\documentclass[review]{elsarticle}
\usepackage{amsfonts,amssymb,amsmath,amsbsy,amsthm,amscd,bbm}
\usepackage[english]{babel}
\usepackage{graphicx}

\newtheorem{theorem}{Theorem}[section]

\newtheorem{corollary}{Corollary}

\def\br{{\bf r}}
\renewcommand{\le}{\leqslant}
\renewcommand{\ge}{\geqslant}

\def\Div{\mbox{div}\,}

\def\bB{{\bf B}}

\def\bE{{\bf E}}

\def\dfrac#1#2{\displaystyle{#1\over #2}}

\def\bv{{\bf V}}
\def\bV{{\bf V}}

\begin{document}

\begin{frontmatter}

\title{Criterion of singularity formation for  radial solutions of the
pressureless Euler-Poisson equations in  exceptional dimension
}


\author{Olga S. Rozanova}
\ead{rozanova@mech.math.msu.su} 

\address[1]{Department of Mechanics and Mathematrics, Moscow
State University, Moscow 119991 Russia}


\begin{abstract}
The spatial dimensions 1 and 4 play an exceptional role for radial solutions of the
pressureless repulsive Euler-Poisson equations. Namely, for any spatial dimension except  1 and 4, any nontrivial solution of the Cauchy problem blows up in a finite time (except in special cases), whereas for  dimensions 1 and 4 there exists a neighborhood of  trivial initial data in the $C^1$ - norm such that the respective solution  preserves the initial smoothness globally. For  dimension 1,  the criterion of the singularity formation in  terms of initial data was known,
 for the case of dimension 4, there was no similar result. In this paper, we close this gap and obtain such a criterion for the case of a more technically complicated case of  dimension 4.
  \end{abstract}


\def\sign{\mathop{\rm sgn}\nolimits}





\end{frontmatter}

\section{Introduction}

In this paper, we study a version of the repulsive Euler-Poisson
equations
\begin{eqnarray}\label{EP}
\dfrac{\partial n }{\partial t} + \Div(n \bv)=0,\quad
\dfrac{\partial \bv }{\partial t} + \left( \bv \cdot \nabla \right)
\bv = \,  \nabla \Phi, \quad \Delta \Phi =n-n_0,
\end{eqnarray}
where the solution components, the scalar functions $n$ (density),
$\Phi$ (force potential), and the vector $\bv$ (velocity) depend on
the time $t$ and the point $x\in {\mathbb R}^d $, $d\ge 1$,
$n_0> 0$ is the density background.

Let $\bE=-\nabla \Phi$.
Under the assumption
 that the solution is sufficiently smooth
  we
 get
\begin{eqnarray*}
n=n_0- \Div \bE,\label{n}\end{eqnarray*} therefore $ n $ can be removed
from the system. Thus, the resulting system is
\begin{eqnarray}\label{4}
\dfrac{\partial \bv }{\partial t} + \left( \bv \cdot \nabla \right)
\bv = \, - \bE,\quad \frac{\partial \bE }{\partial t} + \bv \Div \bE
 = n_0 \bV.
\end{eqnarray}

Denote $ {\bf r}=(x_1, \dots ,x_d), \quad r=|{\bf r}|.$

System \eqref{4}  has an important
class of solutions depending only on
$r$ as
\begin{eqnarray}\label{sol_form}
\bv=F(t,r) \br,\quad \bE=G(t,r) \br,\quad \bB=Q(t,r)\br, \quad
n=n(t,r).
\end{eqnarray}

 Consider the initial data
\begin{equation}\label{CD1}
(\bv, \bE) |_{t=0}= (F_0(r) {\bf r}, G_0(r) {\bf r} ), \quad (F_0(r)
, G_0(r) ) \in C^2(\bar {\mathbb R}_+),
\end{equation}
such that  $n|_{t=0}>0$.

The local in $t$ well-posedness of the the Cauchy problem \eqref{EP}, \eqref{CD1} is discussed, i.e.  in \cite{Tan}. In particular, it is known that the formation of singularity is associated with infinite gradient of the solution.

We call a solution of \eqref{4},  \eqref{CD1} smooth for $t\in
[0,t_*)$, $t_*\le\infty$, if the functions $F$ and $G$ in
\eqref{sol_form} belong to the class $ C^1([0,t_*)\times
\bar{\mathbb R})$. The blow-up of solution implies that the
derivatives of solution tends to infinity as $t\to t_*<\infty$.

The analysis of radially symmetric solutions in the Euler-Poisson equations without pressure and related systems has been the subject of many papers in recent decades; for a review, see \cite{Tan}, \cite{Bhat23}. The cases $n_0=0$ and $n_0>0$ are significantly different, since for $n_0>0$ the motion is oscillatory. The latter case is more difficult and requires a special technique. One of the most interesting and challenging questions is the criterion for the formation of a singularity in terms of initial data. In other words, it is necessary to divide the set of smooth initial data into two parts: solutions with data from one part are globally smooth in time, otherwise the solutions blow up in finite time. It is quite rare that this problem can be completely solved.
For the repulsive case with nonzero background $n_0$ such criterium is known in the case $d=1$
\cite{ELT}, \cite{RChZAMP21}.
It implies that if the initial data
\begin{equation}\label{CD}
(V,E)|_{t=0}=(V_0(x), E_0(x))\in {C^2} ({\mathbb R}).
\end{equation}
are such that for any $x_0\in \mathbb R$ the inequality
\begin{eqnarray} \label {crit2}
\left (V'_0 (x_0) \right) ^ 2 + 2 \, E'_0 (x_0) -n_0 <0
\end {eqnarray}
holds, then the solution of \eqref{4}, \eqref{CD} is periodic in $t$ and  preserves initial smoothness for all $t>0$, otherwise the solution blows up in a finite time. It is easy to see that for $n_0>0$ there exists a neighborhood of the trivial equilibrium $V=E=0$ in the $C^1$-norm such that the solution starting from this neighborhood will be globally smooth.

For the case $d\ge 2$, $d\ne 4$ it is known that, except for special initial data, any nontrivial smooth solution, even a small perturbation of the trivial steady state, blows up in finite time. The proof was done in \cite{R22_Rad} by an analytical method for a small perturbation of the trivial steady state and by a semi-analytical method for an arbitrary case using Floquet theory. Further, in \cite{Carrillo} for a special case and in \cite{Roz_doping} for the general case it was proved that if the Lagrangian trajectories of particles are periodic with different periods, then the corresponding solution necessarily blows up.
In fact, this property is a kind of "folklore" in plasma physics, where \eqref{EP} is one of the most important models.
In other words, the model can have a non-trivial globally smooth solution if and only if the equilibrium equation governing the Lagrangian trajectory of particles has an isochronous center.
For the radially symmetric case, the equation defining the behavior of the trajectories (except for simple waves) is \eqref{defining}, see below, and the Sabatini criterion \cite{Sabatini}, see Appendix, implies that the center $(F=0, \dot F=0)$ is isochronous if and only if $d=1$ and $d=4$.

In this paper, a criterion of the singularity formation   for the remaining case $d=4$ is obtained.
 The possibility of obtaining such a criterion is always based on the existence of first integrals for a nonlinear ODE system.
This is a very non-trivial problem, since there are no general methods for finding such first integrals, even if we suspect their existence in the system. The main difficulty in the case $d=4$ is that the criterion uses not only the values of the derivatives of the data, but also the data themselves.

 It should be noted that the pressureless Euler-Poisson equations with quadratic constraint considered in \cite{Carrillo} are closely related to the \eqref{EP} case, since they can be reduced to \eqref{EP} with $n_0=d$. In this paper, $n_0$ does not depend on the dimension, which is natural for plasma physics \cite{GR75}, \cite{CH18}, we set $n_0=1$. However, all results can be readily reformulated for the case $n_0=d$. The criterion for the formation of a singularity for the case $d= 4$, not considered in \cite{Carrillo}, can be obtained similarly, by the method described here.

The paper is organized as follows. In Sec.\ref{S2} we recall the main technical results of \cite{R22_Rad}. In Sec.\ref{S3} we construct a first integral that helps us to derive the criterion for the singularity formation in  four dimensions in terms of $V_0(r_0), E_0(r_0), {\rm div} V_0(r_0), {\rm div} E_0(r_0) $, $r_0\ge 0$. In Sec.\ref{S4} we derive and analyze the criterion for the important cases of zero initial velocity and zero initial electric field. This halves the dimensionality of the parameter space, since $V_0(r_0)={\rm div} V_0(r_0)=0$ or $E_0(r_0)={\rm div} E_0(r_0)=0$, and the criterion takes a simpler form. We give as an example an initial perturbation of the electric potential in the form of a Gaussian potential (a standard laser pulse) and show that, in contrast to the case of $d=1$, this solution never blows up. We then analyze the profile $E_0(r)$ for which the inequality arising in the criterion turns into an equality. In Sec.\ref{S5} we discuss the possibility of applying the described method to other models.
For the convenience of the readers, we provide in the Appendix the formulations of the theorems that we use in the main text.

\section{Auxiliary results}\label{S2}

This paper is entirely based on the technique of \cite{R22_Rad},
therefore we only list the results that we need.

First of all, we note that from \eqref{sol_form} and \eqref{4}
it follows that $F$ and $G$ satisfy the following Cauchy problem:
 \begin{eqnarray}\label{sys_pol1}
    \dfrac{\partial G}{\partial t}+F r \dfrac{\partial G}{\partial r}=F-4 F
    G, \quad
    \dfrac{\partial F }{\partial t}+F r \dfrac{\partial F}{\partial r}=-F^2 -
    G,
 \end{eqnarray}
\begin{equation*}\label{CD2}
(F(0,r), G(0,r))=(F_0(r), G_0(r)), \quad (F_0(r) , G_0(r) ) \in
C^2(\bar {\mathbb R}_+).
\end{equation*}

 Along the characteristic
 \begin{eqnarray}\label{char}
 \dot r = F r,
 \end{eqnarray}
 starting from the point $r_0\in [0,\infty)$
 system \eqref{sys_pol1} takes the form
  \begin{eqnarray}\label{GF}
\dot G =F-d F
    G,\qquad \dot{F } =  -F^2 -
    G.
\end{eqnarray}

Note that if $r_0=0$, then the corresponding characteristic $r(t)=0$, and all characteristics starting from $r_0>0$ lie in the half-plane $r\ge 0$. Thus, no boundary conditions for $F$ and $G$ at the point $r_0=0$ are required.

We can see that \eqref{GF} can be reduced to
a nonlinear Li\'enard type equation \eqref{Lienard}
\begin{equation}\label{defining}
 \ddot F+(2+d)\,F\, \dot F +F+d \, F^3=0.
 \end{equation}
 The Sabatini criterion \cite{Sabatini} (see Appendix) implies that the oscillations are  isochronous if and only if $d=1$ and $d=4$. The trajectory $r(t)$ can be found from \eqref{char} through $F$ and has the same period as $F$.  Thus, trajectories starting from different points $r_0\in \bar {\mathbb R}_+$ have the same period if and only if $d=1$ or $d=4$, and by \cite{Roz_doping} the trajectories described by \eqref{char}, \eqref{defining} necessarily intersect if $d$ is not equal to 1 or 4.

The case $d=1$ is known \cite{ELT}, \cite{RChZAMP21}, so in what follows we will focus on the case $d=4$.

System \eqref{GF} can be integrated:
\begin{eqnarray}\label{intFG}
&&F^2=\frac{2 G-1}{2}+C_4 |1-4 G|^\frac{1}{2} , \quad
C_4=\frac{1-2 G(0,r_0)+2 F^2(0,r_0)}{2|1-4
G(0,r_0)|^\frac{1}{2}}.
\end{eqnarray}
Further, \eqref{char} and \eqref{GF} imply
\begin{equation}\label{rG}
1-4\,G={\rm const}\, r^{-4}.
\end{equation}
We consider the domain $G< \frac{1}{4}$, since $n>0$.

Along each characteristic, starting from $r_0\in\bar{\mathbb R}_+$, the functions
$F(t)$ and $G(t)$, the solutions of \eqref{GF}, are periodic with period
\begin{equation}\label{T}
  T= 2 \int\limits_{G_-}^{G_+} \frac {d\eta}{(1-4\, \eta )F(\eta)}= 2\pi,
\end{equation}
$F$ is given as   \eqref{intFG}, $G_-<0$ and
$G_+>0$ are the lesser and greater roots of the equation $F(G)=0$.
Moreover, $\int\limits_0^{T} F(\tau) \, d\tau=0$.

 We denote
 \begin{equation}\label{notat}
   {\mathcal D}=\Div \bv,  \quad \lambda=\Div \bE, \quad u={\mathcal D}-4\, F,
\quad v=\lambda-4\, G.
\end{equation}
Evidently, if $u, v, F, G$ are bounded, then
${\mathcal D}$ and $\lambda$ are bounded. Then we get
\begin{equation}\label{uv}
    \dot u=-u^2-2\,F\, u  -v, \quad \dot v=-u\,v + (1-4\, G) \,u- 4\, F\,
    v,
\end{equation}
a quadratically nonlinear system with the coefficient $F$ found from
\eqref{GF}. In fact, \eqref{GF}, \eqref{uv}, can be considered as
a system of 4 ODEs for $G, F, u, v$, where
\eqref{GF} is separated.
If the data \eqref{CD1} are such that the solution to the Cauchy
problem \eqref{uv},
\begin{equation}\label{uv_cd}
u_0(r_0)=({\mathcal D}-4 F)|_{t=0,\,
r=r_0},\quad  v_0(r_0)=({\lambda}-4 G)|_{t=0,\, r=r_0}
\end{equation}
 is
bounded for all fixed $r_0\in [0,+\infty)$, $t\in [0, t_*)$, $t_*\le
\infty$, then ${\mathcal D}$ and $\lambda$ are bounded and for $t\in
[0, t_*)$  there exists a $C^1$-smooth  solution to \eqref{4},
\eqref{CD1} with a positive density.

 According to  the Radon lemma
\cite{Riccati}, Theorem 3.1, \cite{Radon}  (see Appendix), system \eqref{uv} can be linearized, so
we obtain a linear Cauchy problem
\begin{eqnarray}
\label{matr}
 \begin{pmatrix}
  \dot q\\
  \dot p_1\\
  \dot p_2\\
\end{pmatrix}
=\begin{pmatrix}
0& 1& 0\\
0&- 2\,F& -1\\
0& 1-4\,G  & -4\,F\\
\end{pmatrix}
\begin{pmatrix}
  q\\
  p_1\\
  p_2\\
\end{pmatrix},\quad
\begin{pmatrix}
  q\\
  p_1\\
  p_2
  \\
\end{pmatrix}(0)=\begin{pmatrix}
  1\\
  u_0\\
  v_0\\
\end{pmatrix},
\end{eqnarray}
with periodical coefficients, known from \eqref{GF}. System
\eqref{matr} implies 
\begin{eqnarray*}\label{p1}
  \ddot p_1+6 F \dot p_1 +(1+6 F^2 -6 G) p_1=0.
\end{eqnarray*}
The standard change of the variable $p_1(t)=P(t)\, e^{-
4\int\limits^t_0 F(\tau)\,d\tau}$ reduces the latter
equation to
\begin{eqnarray}\label{P}
  \ddot P+{Q} P=0, \quad {Q}=1-3 G.
\end{eqnarray}

 The solution of \eqref{uv} blows up if
and only if $q(t)$ vanishes at some point $t_*,\,0<t_*<\infty$.

From \eqref{matr} we find
\begin{equation}\label{q}
 q(t)=1+\int\limits_0^t p_1(\tau)\,d
\tau=1+\int\limits_0^t P(\xi)\, e^{- 3\int\limits^\xi_0
F(\tau)\,d\tau}\,d\xi.
\end{equation}

\section{First integral and the main theorem}\label{S3}

As in the $d=1$ case, the key point in this new situation is the existence of a first integral relating the derivatives of the solution.
To simplify the notation, we set
\begin{equation}\label{denot1}
M=\sqrt{1-4 G}, \quad S= F M^2.
\end{equation}
Since $G<\frac14$, then $M>0$.

Due to \eqref{intFG} and \eqref{rG},
\begin{equation}\label{denot2}
 S= \pm \frac{ M^2}{2} \sqrt{4 C_4 M-1-M^2}, \quad r=\frac{K}{\sqrt{M}},\quad Q= \frac{1+M^2}{4},
\end{equation}
with
\begin{equation}\label{const}
 C_4= \frac{1+M_0^2+4 F_0^2}{4 M_0}\ge \frac12, \quad K=r_0 \sqrt{M_0}\ge 0, \quad M_0=\sqrt{1-4G_0}\le 1.
\end{equation}

Note that if $G\in [G_-, G_+]$ (see \eqref{T}), then $0<M_-\le M\le M_+$, where
\begin{equation}\label{Mpm}
M_\pm= 2 C_4 \pm \sqrt{4 C_4^2-1}=\frac{1+M_0+4 F_0^2 \pm \sqrt{(1+M_0+4 F_0^2)^2- 4 M_0^2}}{2 M_0}.
\end{equation}

If we denote $R(t)=\dot P(t), $ then we can write \eqref{GF}, \eqref{P}, \eqref{q} (using \eqref{char}) as
\begin{equation*}\label{PRGq}
\dot P=R,\quad \dot R=-Q P,\quad \dot G=S, \quad \dot q =\frac{ r_0^3}{r^3}\,P,
\end{equation*}
or (see \eqref{denot1}, \eqref{denot2})
\begin{eqnarray}\label{PRq}
&& P'(M)=-\frac{M R(M)}{2 S(M)}, \nonumber\\ && R'(M)  =\frac{M P(M) Q(M)}{2 S(M)},\label{PRq}\\ && q'(M)= -\frac{M P(M) r_0^3}{2  r^3 (M) S(M)}.\nonumber
\end{eqnarray}
The main point that ensures the success of obtaining the criterion is that
system \eqref{PRq} can be integrated:
\begin{eqnarray}
q(M) &=& C_1 + \frac{1}{K_1^3}\,   \left[  C_2\,{\rm sign} \, F{\sqrt{4 C_4 M-1-M^2}}+{C_3}{M}\right],\label{YM}\\
P(M)& =& -\frac {C_2 (2 C_4 -M) - {\rm sign} \, F\, C_3 {\sqrt{4 C_4 M-1-M^2}}}{\sqrt{M}} \label{PM} \\
R(M) & =&  -\frac{ {\rm sign} \, F\,C_2 (2 C_4 +M)  {\sqrt{4 C_4 M-1-M^2}} + C_3 (1-M^2)  }{2\sqrt{M}}, \label{RM}
\end{eqnarray}
where $K_1=\sqrt{M_0}$, constants $C_1$, $C_2$ and $C_3$ depend on $q(M_0)$, $P(M_0)$, $R(M_0)$, i.e. on $G_0$, $F_0$,  $u_0$, $v_0$.
The motion is periodic in $t$. If $F_0>0$, then $M$ increases from $M_0$ to $M_+$ (${\rm sign} F =1$), then at the point $M_+$ the sign of $F$ changes and $M$ decreases from $M_+$ to $M_-$.

The constants $C_2$ and $C_3$ can be found from the linear inhomogeneous system \eqref{PM}, \eqref{RM} if instead of $M$, $P(M)$, $R(M)$ we substitute $M_0$, $P_0=P(M_0)$, $R_0=R(M_0)$, respectively. The determinant is
\begin{eqnarray}\label{deter}
D=\frac{1}{M_0}\left((1+4 M_0+M_0^2) C_4 -(M_0^2+M_0+1)    \right)\ge \frac{(M_0^2-1)^2}{4 M_0^2}\ge 0.
\end{eqnarray}
We used the fact that $C_4\ge \frac{1+M_0^2}{4 M_0}$. The point $M_0=1$ corresponds to $G_0=F_0=0$, so $D>0$ for all non-trivial solutions.

Let us denote
\begin{eqnarray*}\nonumber
Y(M)=\frac{1}{K_1^3}\,   \left[  C_2\,{\rm sign} \, F{\sqrt{4 C_4 M-1-M^2}}+{C_3}{M}\right].\label{Y}
\end{eqnarray*}
Computations show that
\begin{eqnarray}\label{C1}
&& C_1= 1- Y(M_0),
\end{eqnarray}
where
\begin{eqnarray}\label{C2}
&& C_2=  - 2 M_0^{\frac32}\frac{4 F_0 R_0+(M_0^2-1) P_0 }{((M_0+1)^2+4 F_0^4)((M_0-1)^2+4 F_0^4)}\\
&& C_3=   2 M_0^{\frac12}\frac{(4 F_0^2 + 3 M_0^2 +1) F_0 P_0+(M_0^2-4 F_0^2-1) R_0 }{((M_0+1)^2+4 F_0^4)((M_0-1)^2+4 F_0^4)}.
 \label{C3}
\end{eqnarray}

It can be readily shown that
\begin{eqnarray}\label{P0}
&& P(0)=P(M_0)=p_1(0)=u_0, \\ && R(0)=R(M_0)=-p_2(0)+2 F(0) p_1(0)=-v_0 +2 F_0 u_0. \label{R0}
\end{eqnarray}

 To find solutions that blow up, we need to find a set of initial data such that $q(M)\le 0$ for some $M\in [M_-,M_+]$, subject to $q(M_0)=1$.

The function $q(M)$ (see \eqref{YM}) reaches its extrema at the points where $q'(M)$, i.e. where $P(M)=0$ (see \eqref{PRq}). Thus, as follows from \eqref{PM}, these points can be found as a solution of the equation
\begin{equation*}\label{M_ext}
 C_2 (2 C_4 -M) -C_3 \, {\rm sign} F \, {\sqrt{4 C_4 M-1-M^2}}=0.
\end{equation*}
This equation has a solution,
\begin{equation}\label{MM}
 M^*_\pm= 2 C_4 - {\rm sign} F \, \lambda \,\sqrt{4 C_4^2-1}, \quad \lambda =\frac{C_3}{\sqrt{C_2^2+C_3^2}}, \quad M_-\le M^*_\pm\le M_+.
\end{equation}
As follows from \eqref{Mpm}, for $C_2>0$ we have $P(M_-)<0$, $P(M_+)>0$, therefore for $F>0$ the function $q(M)$ has a minimum at the point $M^*_+$. Otherwise, for $C_2<0$ the minimum is at the point  $M^*_-$, where $F<0$. If $C_2=0$, then the extrema are at the boundary points of $M_\pm$.

Note that due to the periodicity of the function $q(t)$ (see \eqref{q}), if it vanishes at some point, then this point lies inside its period. Recall that the period is equal to $2\pi$ (see \eqref{T}).

Thus, we obtain the following theorem.
\begin{theorem}\label{T1}
Let the data \eqref{CD1} are such that for all $r_0\ge 0$
\begin{eqnarray}\label{crit}
q^*=\min\{C_1-Y(M_+^*), \,C_1-Y(M_-^*) \}>0.
\end{eqnarray}
Then the solution of problem \eqref{4}, \eqref{CD1} preserves smoothness for all $t>0$. Otherwise, the solution blows up  in finite time $0<T_*< 2 \pi$.

All constants in \eqref{crit} are presented as \eqref{C1}, \eqref{C2}, \eqref{C3}, \eqref{MM}, \eqref{P0}, \eqref{R0}, \eqref{uv_cd}.
\end{theorem}

{\it Remark 1.} Since $C_1$ tends to 1, and $C_2$, $C_3$ tend to zero as $u_0$, $v_0$ (see \eqref{C1}, \eqref{C2}, \eqref{C3}, \eqref{P0}, \eqref{R0})
we see from \eqref{crit} that there exists a neighborhood of the trivial equilibrium $V=E=0$ in the $C^1$-norm such that the solution starting from this neighborhood is smooth for all $t>0$.

\medskip

Despite the explicit algorithm, it is not possible to conveniently represent the set of points corresponding to a globally smooth solution in the 4-dimensional space $(G_0, F_0, u_0, v_0)$.

However, for important special cases, simpler results can be obtained.

\section{Zero initial velocity}\label{S4}

\subsection{Zero initial velocity}\label{S41}

\begin{corollary} Let the initial data \eqref{CD1} be such that $V_0=0$. Then the solution \eqref{4}, \eqref{CD1} preserves global smoothness in $t$ if and only if for each point $r_0\ge 0$ the following property holds:
\begin{equation}\label{critV0}
  {\rm div} E_0(r_0) > 6 \, G_0(r_0)- \frac12.
\end{equation}
Here $G_0(r_0)=\frac{E_0(r_0)}{r_0}<\infty$.
\end{corollary}
\proof
If $V_0\equiv0$, then $F_0=u_0=0$, therefore according to \eqref{C1} - \eqref{R0} we have $C_1=1+\frac {2 v_0 }{(M_0^2-1)}$, $C_2=0$, $C_3=-\frac{2 v_0 \sqrt{M_0}}{M_0^2-1}$. Here $M_0\ne 1$, see the remark after \eqref{deter}.
Further, according to \eqref{Mpm} we have $M_\pm=M_0, \, \frac{1}{M_0}$, and $M^*_\pm=M_\pm$, see the remark after \eqref{MM}.
Let us find $q^*$. It is easy to check that $C_1-Y(M_0)=1$ and $C_1-Y\left(\frac{1}{M_0}\right)=1+\frac{2 v_0}{M_0^2}$,
so the criterion looks like
 \begin{equation*}\label{critV0M}
 v_0> -\frac12 \,M_0^2,
\end{equation*}
which coincides with \eqref{critV0} taking into account \eqref{notat} and \eqref{const}. $\Box$


\medskip

{\it Remark 2.}
Note that in the case $d=1$, as follows from \eqref{crit2} for $V_0=0$, $n_0=1$, the analogous domain on the plane $({\rm div}E_0, G_0)$ is ${\rm div}E_0=E'_0<\frac12$, it does not depend on $G_0<1$. Thus, in contrast to the case $d=4$, the smoothness region containing the origin is bounded on the other side.

\medskip

{\it Example 1.} For  example, we choose the initial data in the form of a standard laser pulse
 \cite{esarey09}, \cite{CH18}:
\begin{equation}\label{DE}
E_0= a \, {\bf r} \, e^{-\frac{r^2}{2}}, \qquad V_0=0.
\end{equation}
Here $G_0(r)=a  \, e^{-\frac{r^2}{2}}$, $a<\frac14$ (to guarantee positivity of the initial density).

It can be readily computed that $q^*(r)=\frac12- (r^2+2)\,a  \, e^{-\frac{r^2}{2}}>0$ for all possible values of $a$, in other words, the standard laser pulse never blows up.

\medskip

{\it Remark 2.} Let us find the critical form $G_0^*(r)$ of the function $G_0(r)$, that is, the function that is the solution of the equation $q^*=0$,
\begin{equation*}\label{Gcrit}
r \dfrac{dG^*_0(r)}{dr}= 2 G_0^*(r) - \frac12.
\end{equation*}
The solution is $G^*_0(r)=\frac14- C\,r^2$, $C={\rm const}>0$.


\medskip

\subsection{Zero initial electric field}\label{S42}

\begin{corollary} Let the initial data \eqref{CD1} be such that $E_0=0$. Then the solution \eqref{4}, \eqref{CD1} preserves global smoothness in $t$ if and only if for each point $r_0\ge 0$ the following property holds:
\begin{equation}\label{critE0}
 \min \left[1 -\left(F_0(r_0) \pm \frac{F_0^2(r_0)-1}{\sqrt{F_0^2(r_0)+1}}\right)\, \left({\rm div} V_0(r_0) -4 F_0(r_0)\right) \right] >0.
\end{equation}
Here $F_0=\frac{V_0(r_0)}{r_0}<\infty$.
\end{corollary}
\proof
 In this case $M_0=1,$ $u_0=1$, $ M^*_\pm=1+2 F_0^2 \pm 2 F_0 \, \sqrt{1+F_0^2}$, $C_1=1-\frac12 \,\frac{1+3 F_0^2}{1+F_0^2} \,u_0$,
 $C_2=-\frac{4 F_0^2 u_0}{(1+2 F_0^2)^2-2}$,  $C_3=\frac{2 F_0 (1-F_0^2) u_0}{(1+2 F_0^2)^2-2}$. Computations show that
 \begin{equation*}
 q^* = \min \left[1 -\left(F_0 \pm \frac{F_0^2-1}{\sqrt{F_0^2+1}}\right)\, u_0 \right],
\end{equation*}
which implies \eqref{critE0}. $\Box$

\begin{figure}[htb]
\begin{minipage}{0.5\columnwidth}
\includegraphics[scale=0.35]{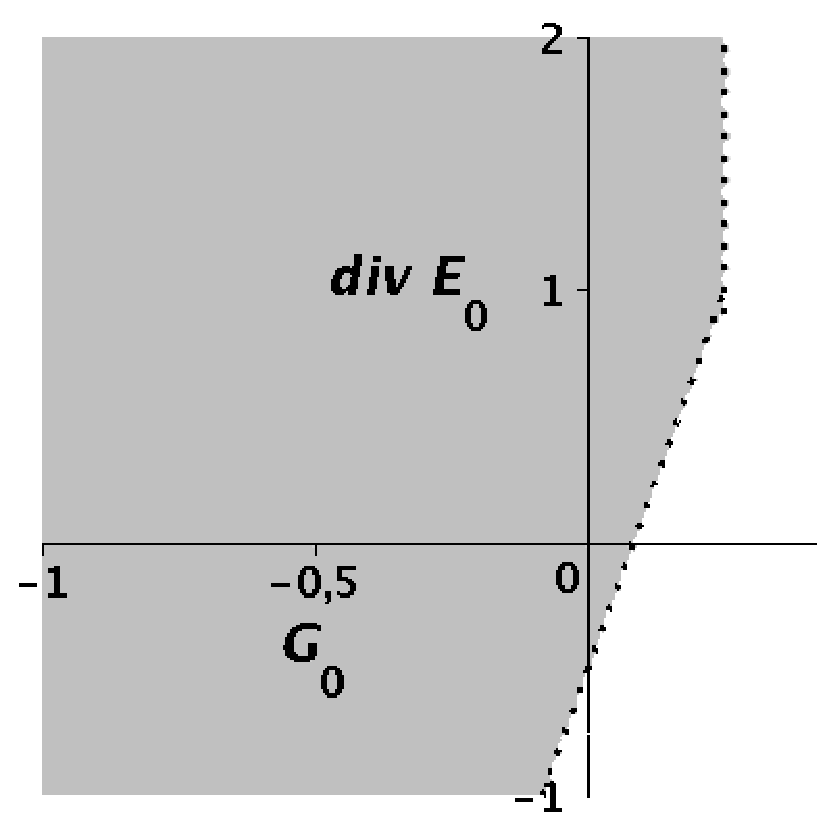}
\end{minipage}
\hspace{0.7cm}
\begin{minipage}{0.5\columnwidth}
\includegraphics[scale=0.35]{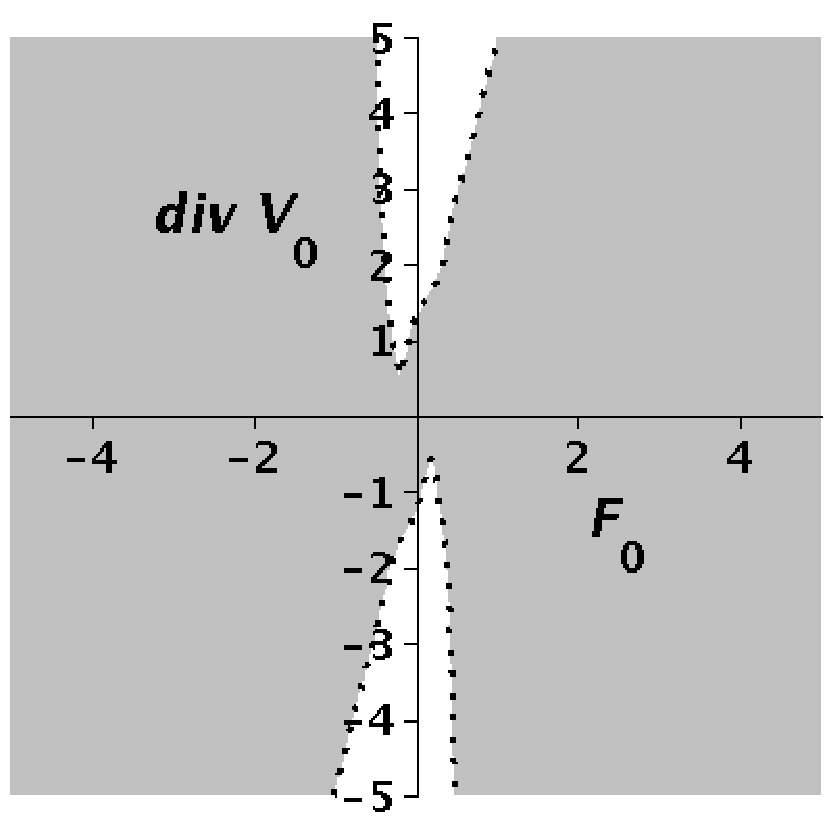}
\end{minipage}
\caption{The domains of smoothness (shaded) for the case of zero initial velocity on the plane ($(G_0, {\rm div} E_0)$) (left) and for the case of zero initial electric field on the plane ($(F_0, {\rm div} V_0)$) (right).}\label{Pic1}
\end{figure}

Fig.1 presents the domains of smoothness for the case $V_0=0$ on the plane ($ G_0, {\rm div} E_0$), defined by inequalities \eqref{critV0} and $G_0<\frac14$   (left), and for the case $E_0=0$ on the plane ($ F_0, {\rm div} V_0$), defined by inequality \eqref{critE0} (right). These figures, in particular, give an idea that the domains of smoothness in the full four-dimensional space of parameters ($ G_0, F_0,  {\rm div} E_0, {\rm div} V_0$)   is structured in a  complex way.

\section{Discussion}\label{S5}

Let us describe other models that have  radially symmetric solutions and can be considered using the methods presented above.

1. First of all, it is  the pressureless Euler-Poisson equations with quadratic confinement \cite{Carrillo}, which has the form
\begin{eqnarray}\label{CC}
\dfrac{\partial n }{\partial t} + \Div(n \bv)=0,\quad
\dfrac{\partial \bv }{\partial t} + \left( \bv \cdot \nabla \right)
\bv = -\int_{{\mathbb R}^d}\, \nabla_{\bf x} N({\bf x}-{\bf y}) n(t,{\bf y}) \, d{\bf y} -{\bf x},
\end{eqnarray}
where $-N({\bf x}) $ is the fundamental solution of the Laplace operator in $d$ - dimensional space, i.e. $-\Delta N ({\bf x})=\delta({\bf x})$. Here, as before, $n(t,{\bf x})$ and $\bv (t,{\bf x})$ are  the density and the velocity, $d\ge 2$. It can be shown that the right hand side term in the second equation \eqref{CC} is $-\nabla \Psi$, where $\Psi $ is the solution of   $\Delta \Psi =\rho-d$. Thus, it coincide with \eqref{EP} for $n_0=d$.

We can make the same computations that in \cite{R22_Rad} with the change $n_0=1$ to $n_0=d$ and obtain
along the characteristic \eqref{char}
  $r_0\in [0,\infty)$
 system
  \begin{eqnarray*}\label{GF1}
\dot G =dF-d F
    G,\qquad \dot{F } =  -F^2 -
    G.
\end{eqnarray*}
Thus,
instead of  \eqref{defining} we get the Li\'enard type  equation \eqref{Lienard} as
\begin{equation*}\label{defining1}
 \ddot F+(2+d)\,F\,\dot  F +d F+d \, F^3=0,
 \end{equation*}
for which the Sabatini criterion yields that the equilibrium $(F=0,\dot F=0)$ is an isochronous center if and only if $d=1$ and $d=4$.
If the solution is a simple wave (i.e. $F=F(G)$), then the situation is completely different. Instead of \eqref{sys_pol1} we get
\begin{eqnarray*}\label{F(G)1}
    \dfrac{\partial G}{\partial t}+F(G) r \dfrac{\partial G}{\partial r}=dF(G)(1- G),
    \end{eqnarray*}
    and
   the periods of oscillations
  are equal for all characteristics (see Sec.6 of \cite{R22_Rad}).
This means that the globally smooth solutions obtained in \cite{Carrillo} for $d\ge 2$, $d\ne 4$ are simple waves.

\medskip

2. Another model, which is related to the pressureless Euler-Poisson equations, is a system with a quadratic confinement of a more general form than \eqref{CC},
\begin{eqnarray*}\label{EP2}
\dfrac{\partial n }{\partial t} + \Div(n \bv)=0,\quad
\dfrac{\partial \bv }{\partial t} + \left( \bv \cdot \nabla \right)
\bv = \,  -\nabla \Phi- k {\bf x}, \quad \Delta \Phi =n-n_0,
\end{eqnarray*}
where $k$ is some constant. On the right-hand side of the equation for $\bv$, one can also write a non-local term similar to \eqref{CC}.
For $n_0> k d$\, the motion is oscillatory, and $d=1$, $d=4$ are also critical dimensions, which allows the existence of global in time non-trivial solutions with data from the neighborhood of the trivial steady state. For $n_0\le k d$ the solutions are not oscillatory, but their behavior can still be described using a technique based on the Radon lemma.

\medskip

3. Note that if there is a quantity $S$ that is conserved along Lagrangian trajectories,
\begin{eqnarray*}\label{EP2}
\dfrac{\partial S }{\partial t} +  \bv\cdot \nabla S=0,
\end{eqnarray*}
then together with the continuity equation we have an additional conservation law
\begin{eqnarray*}
\frac{\partial n S }{\partial t} + \Div(n S \bv)=0.
\end{eqnarray*}
If we change the elliptic part of the Euler-Poisson system to
\begin{eqnarray*}
\Delta \Phi =n S-n_0 S_0,\quad n_0 S_0= {\rm const},
\end{eqnarray*}
then, by introducing $\bE=-\nabla \Psi$, we can again obtain a system like \eqref{4} and apply the same methods as before.

In particular, any reasonable function $\Phi(|{\bf M}|)$, ${\bf M}=(M_1,\dots, M_d)$, $M_i=\int\limits_0^{x_i} \,n(t,{\bf y}) dy_i$, ${\bf y}\in {\mathbb R}^d$, $i=1,\dots, d$, can be chosen as $S$.

\medskip

\section*{Acknowledgements}

 Supported by Russian Science Foundation  grant 23-11-00056 through RUDN University.

\section*{Appendix}

1. The Sabatini criterion (1999).
\begin{theorem}\label{S}
Let us consider a  Li\'enard type equation
\begin{equation}\label{Lienard}
\ddot y+ f (y) \dot y+g(y)=0,
\end{equation}
where $f, g$  are analytic, $g$ odd, $f (0)=g(0)=0$, $g'(0)>0$. Then
$\mathcal O =(y,\dot y)=(0,0)$ is a center if and only if $f$ is odd and
$\mathcal O )$ is an isochronous center if and only if
\begin{equation}\label{tau}
\tau(y) :=\left(\int\limits_0^y sf(s) ds \right)^2-y^3 (g(y)-g'(0)y)=0.
\end{equation}
\end{theorem}

\medskip

2. The Radon lemma (1927).
\begin{theorem}
\label{T2} A matrix Riccati equation
\begin{equation}
\label{Ric}
 \dot W =M_{21}(t) +M_{22}(t)  W - W M_{11}(t) - W M_{12}(t) W,
\end{equation}
 {\rm (}$W=W(t)$ is a matrix $(n\times m)$, $M_{21}$ is a matrix $(n\times m)$, $M_{22}$ is a matrix  $(m\times m)$, $M_{11}$ is a matrix  $(n\times n)$, $M_{12} $ is a matrix $(m\times n)${\rm )} is equivalent to the homogeneous linear matrix equation
\begin{equation}
\label{Lin}
 \dot {\mathcal Y} =M(t) {\mathcal Y}, \quad M=\left(\begin{array}{cc}M_{11}
 & M_{12}\\ M_{21}
 & M_{22}
  \end{array}\right),
\end{equation}
 {\rm (}${\mathcal Y}={\mathcal Y}(t)$  is a matrix $(n\times (n+m))$, $M$ is a matrix $((n+m)\times (n+m))$ {\rm )} in the following sense.

Let on some interval ${\mathcal J} \in \mathbb R$ the
matrix-function $\,{\mathcal Y}(t)=\left(\begin{array}{c}\mathfrak{Q}(t)\\ \mathfrak{P}(t)
  \end{array}\right)$ {\rm (}$\mathfrak{Q}$  is a matrix $(n\times n)$, $\mathfrak{P}$  is a matrix $(n\times m)${\rm ) } be a solution of \eqref{Lin}
  with the initial data
  \begin{equation*}\label{LinID}
  {\mathcal Y}(0)=\left(\begin{array}{c}I\\ W_0
  \end{array}\right)
  \end{equation*}
   {\rm (}$ I $ is the identity matrix $(n\times n)$, $W_0$ is a constant matrix $(n\times m)${\rm ) } and  $\det \mathfrak{Q}\ne 0$ on ${\mathcal J}$.
  Then
{\bf $ W(t)=\mathfrak{P}(t) \mathfrak{Q}^{-1}(t)$} is the solution of \eqref{Ric} with
$W(0)=W_0$ on ${\mathcal J}$.
\end{theorem}

\end{document}